\documentclass[12pt]{elsart}
\usepackage{amsfonts}
\usepackage{amsmath}
\usepackage{epic}
\usepackage{eepic}

\advance \topmargin by -\headheight
\advance \topmargin by -\headsep
\evensidemargin \oddsidemargin
\DeclareMathOperator{\Ad}{Ad}

\DeclareMathOperator{\cyl}{cyl}
\DeclareMathOperator{\diag}{diag}
\DeclareMathOperator{\dist}{dist}
\DeclareMathOperator{\rk}{rank}

\DeclareMathOperator{\spec}{sp}

\newcommand{\nats}{{\mathbb N}}
\newcommand{\ints}{{\mathbb Z}}
 
\newcommand{\cplx}{{\mathbb C}}
\newcommand{\rats}{{\mathbb Q}}
\newcommand{\torus}{{\mathbb T}}

\newcommand{\Edge}{{\mathcal E}}

\newcommand{\iver}{{\tt i}}
\newcommand{\tver}{{\tt t}}
\newcommand{\erange}{{\mathcal R}}

\newtheorem{theorem}{Theorem}[section]
\newtheorem{lemma}[theorem]{Lemma}

\newtheorem{proposition}[theorem]{Proposition}
\newtheorem{remark}[theorem]{Remark}
\newtheorem{example}[theorem]{Example}

\newtheorem{definition}[theorem]{Definition}

\newenvironment{proof}{{\em Proof. }}{\qed}

\begin{document}
\begin{frontmatter}
\title{The Rohlin Property for Shifts of Finite Type}
\author{Charles G. Holton}
\ead{cholton@math.utexas.edu}
\address{Department of Mathematics\\ University of Texas at Austin\\
1 University Station C1200\\
Austin TX 78712-0257, US}

\begin{abstract}
We show that an automorphism of a unital AF $C^*$-algebra with the
approximate Rohlin property has the Rohlin property.  This generalizes a
result of Kishimoto.  Using this we show that the shift automorphism on 
the bilateral $C^*$-algebra associated with an aperiodic irreducible 
shift of finite type has the Rohlin property.
\end{abstract}

\end{frontmatter}
\section{Introduction}

The Rohlin Lemma is one of the foundations of ergodic theory.  The simple 
idea of levels of a stack moving almost cyclically under a transformation has 
proved to be an important and powerful tool.  In 1975, A. Connes proved a 
``noncommutative analogue of the very useful tower theorem of Rohlin'' 
for automorphisms of von Neumann algebras (\cite{Co}).  The notion of a 
Rohlin property for automorphisms of $C^*$-algebras is due to Herman and 
Ocneanu (\cite{HO}).  It is useful in the study of automorphisms up 
to inner automorphism (see, for example, 
\cite{BEK2,BKRS,EK,Ki0,Ki1,Ki2,Ki3}).  Evans and Kishimoto 
showed in \cite{EK} that two 
automorphisms of an AF $C^*$-algebra having the Rohlin property are outer 
conjugate if they induce the same action on $K_0.$
At present, only a few classes of automorphisms of 
$C^*$-algebras are known to possess the Rohlin property.

Our work deals with the bilateral $C^*$-algebras associated with shifts of 
finite type, studied in \cite{Kr,Ru,Wa1}.
The shift map on the shift space induces an automorphism of the $C^*$-algebra.
Our main result is that these automorphisms have the Rohlin property. 
It
follows from the result of Evans and Kishimoto in \cite{EK} that the crossed-product of the bilateral $C^*$-algebra with the $\ints$-action of the 
automorphism is determined up to isomorphism by its $K$-theory.  
Krieger showed in
\cite{Kr} that $K$-theory for shifts of finite type is an invariant 
for shift equivalence, and thus crossed-product algebras arising from shift 
equivalent shifts of finite type are isomorphic.

We define the Rohlin property and the approximate Rohlin property in 
Section~\ref{sec-appx}.  We also prove a few lemmas about projections and
partial isometries in $C^*$-algebras for use in the later sections.  
In Section~\ref{sec-gen} we prove that for unital AF algebras the 
approximate Rohlin property implies the Rohlin property.  This improves the
result of Kishimoto in \cite{Ki1}, where it is shown for for unital simple 
AF algebras whose $K_0$ group has finite rank and no infinitesimals.  We need 
this more general result because the $K_0$ groups for the algebras we are 
considering usually have infinitesimals.  We show that even less than the
approximate Rohlin property is needed to deduce the Rohlin property.

The construction and basic properties of the bilateral $C^*$-algebra 
and the automorphism associated with a shift of finite type are given 
in Section~\ref{sec-def}.  This is the algebra associated
with one of Krieger's $K$-theoretic invariants in \cite{Kr}, although the
algebra itself is not discussed there.  It is a special case of Ruelle's 
construction in \cite{Ru}.  Wagoner studied this and two other algebras 
in \cite{Wa1}.

We discuss shift equivalence and the
relationship with the automorphism in Section~\ref{sec-SE}.  
We show that $K_0$ for a bilateral $C^*$-algebra has no 
infinitesimals if and only if the shift of finite type is shift equivalent 
to a full shift.  Finally, in Section~\ref{sec-rp} we prove that the 
automorphism of the 
bilateral $C^*$-algebra associated with an irreducible aperiodic shift 
of finite type has the Rohlin property.  (This may be deduced for full shifts
from the result of Kishimoto in \cite{Ki4}.)  It follows that the crossed 
product
algebra is no stronger an invariant than the $K$-theoretic invariants introduced
in \cite{Kr}.  It is in fact weaker because one can easily show that the full 
$2$ shift and the full $4$ shift have isomorphic crossed product algebras.
We still do not know whether the pair consisting of the bilateral $C^*$-algebra
and the automorphism associated with a shift of finite type is a stronger 
invariant than shift equivalence.

\section{Definitions and basic lemmas}\label{sec-appx}
Let $A$ be a unital $C^*$-algebra and write $1$ for the identity of $A.$  
We quote a definition from \cite{Ki1}.

\begin{definition}\label{def-rp}An automorphism $\alpha$ of a unital $C^*$-algebra $A$ has the Rohlin property if for any $k\in\nats$ there are positive integers $k_1,\dots,k_m\geq k$ satisfying the following condition: For any finite subset $F$ of $A$ and $\varepsilon>0$ there are projections $e_{i,j},$ $i=0,\dots,m-1,$ $j=0,\dots,k_i-1$ in $A$ such that
\begin{enumerate}
\item $\displaystyle\sum_{i=0}^{m-1}\sum_{j=0}^{k_i-1}e_{i,j}=1,$
\item $\|\alpha(e_{i,j})-e_{i,j+1}\|<\varepsilon,$
\item $\|[x,e_{i,j}]\|<\varepsilon,$
\end{enumerate}
for $i=0,\dots,m-1,$ $j=0,\dots,k_i-1,$ and $x\in F,$ where $e_{i,k_i}=e_{i,0}.$
\end{definition}

\noindent{\bf Notation. }A {\em stack of height $m$} for $\alpha$ is a family of orthogonal projections $f_0,f_1,\dots,f_{m-1}$ such that $\alpha(f_i)=f_{i+1},$ $i=0,1,\dots,m-2.$  By an $\varepsilon$-{\em approximate stack} we mean a family of orthogonal projections $f_0,f_1,\dots,f_{m-1}$ such that $\|\alpha(f_i)-f_{i+1}\|<\varepsilon,$ $i=0,1,\dots,m-2.$  We call a stack [or an $\varepsilon$-approximate stack] $f_0,f_1,\dots,f_{m-1}$ {\em cyclic} if $\alpha(f_{m-1})=f_0$ [resp. if $\|\alpha(f_{m-1})-f_0\|<\varepsilon$].

\begin{figure}
\setlength{\unitlength}{1mm} {
\renewcommand{\dashlinestretch}{2}
\renewcommand{\drawlinestretch}{1000}
\begin{center}
\begin{picture}(125,55)(0,-6)
\drawline(0,0)(20,0)
\drawline(30,0)(50,0)
\drawline(80,0)(100,0)
\drawline(0,-2)(0,2)
\drawline(20,-2)(20,2)
\drawline(30,-2)(30,2)
\drawline(50,-2)(50,2) 
\drawline(80,-2)(80,2)
\drawline(100,-2)(100,2)

\put(10,-1){\makebox(0,0)[ct]{$f_0$}}
\put(40,-1){\makebox(0,0)[ct]{$f_m$}}
\put(90,-1){\makebox(0,0)[ct]{$f_{(n-1)m}$}}

\put(65,0){\makebox(0,0)[ct]{$\dots$}}

\put(10,7.5){\makebox(0,0)[ct]{$\uparrow$}}
\put(12.5,6.5){\makebox(0,0)[ct]{$\alpha$}}
\put(40,7.5){\makebox(0,0)[ct]{$\uparrow$}}
\put(42.5,6.5){\makebox(0,0)[ct]{$\alpha$}}
\put(90,7.5){\makebox(0,0)[ct]{$\uparrow$}}
\put(92.5,6.5){\makebox(0,0)[ct]{$\alpha$}}

\drawline(0,15)(20,15)
\drawline(30,15)(50,15)
\drawline(80,15)(100,15)
\drawline(0,13)(0,17)
\drawline(20,13)(20,17)
\drawline(30,13)(30,17)
\drawline(50,13)(50,17) 
\drawline(80,13)(80,17)
\drawline(100,13)(100,17)

\put(10,14){\makebox(0,0)[ct]{$f_1$}}
\put(40,14){\makebox(0,0)[ct]{$f_{m+1}$}}
\put(90,14){\makebox(0,0)[ct]{$f_{(n-1)m+1}$}}

\put(65,15){\makebox(0,0)[ct]{$\dots$}}

\put(10,22.5){\makebox(0,0)[ct]{$\uparrow$}}
\put(12.5,21.5){\makebox(0,0)[ct]{$\alpha$}}
\put(40,22.5){\makebox(0,0)[ct]{$\uparrow$}}
\put(42.5,21.5){\makebox(0,0)[ct]{$\alpha$}}
\put(90,22.5){\makebox(0,0)[ct]{$\uparrow$}}
\put(92.5,21.5){\makebox(0,0)[ct]{$\alpha$}}

\put(10,30){\makebox(0,0)[ct]{$\vdots$}}
\put(40,30){\makebox(0,0)[ct]{$\vdots$}}
\put(90,30){\makebox(0,0)[ct]{$\vdots$}}

\put(10,35){\makebox(0,0)[ct]{$\uparrow$}}
\put(12.5,34){\makebox(0,0)[ct]{$\alpha$}}
\put(40,35){\makebox(0,0)[ct]{$\uparrow$}}
\put(42.5,34){\makebox(0,0)[ct]{$\alpha$}}
\put(90,35){\makebox(0,0)[ct]{$\uparrow$}}
\put(92.5,34){\makebox(0,0)[ct]{$\alpha$}}

\drawline(0,45)(20,45)
\drawline(30,45)(50,45)
\drawline(80,45)(100,45)

\drawline(0,43)(0,47)
\drawline(20,43)(20,47)
\drawline(30,43)(30,47) 
\drawline(50,43)(50,47)
\drawline(80,43)(80,47)
\drawline(100,43)(100,47)

\put(10,44){\makebox(0,0)[ct]{$f_{m-1}$}}
\put(40,44){\makebox(0,0)[ct]{$f_{2m-1}$}}
\put(90,44){\makebox(0,0)[ct]{$f_{nm-1}$}}

\put(65,45){\makebox(0,0)[ct]{$\dots$}}

\put(115,-1){\makebox(0,0)[ct]{$e_0$}}
\put(115,14){\makebox(0,0)[ct]{$e_1$}}
\put(115,44){\makebox(0,0)[ct]{$e_{m-1}$}}
\put(115,7.5){\makebox(0,0)[ct]{$\uparrow$}}
\put(117.5,6.5){\makebox(0,0)[ct]{$\alpha$}}
\put(115,22){\makebox(0,0)[ct]{$\uparrow$}}
\put(117.5,21.5){\makebox(0,0)[ct]{$\alpha$}}
\put(115,35){\makebox(0,0)[ct]{$\uparrow$}}
\put(117.5,34){\makebox(0,0)[ct]{$\alpha$}}
\put(115,30){\makebox(0,0)[ct]{$\vdots$}}

\end{picture}
\end{center}
}
\caption{The stack from Remark~\ref{rem-collapse}.}
\label{fig-collapse}
\end{figure}
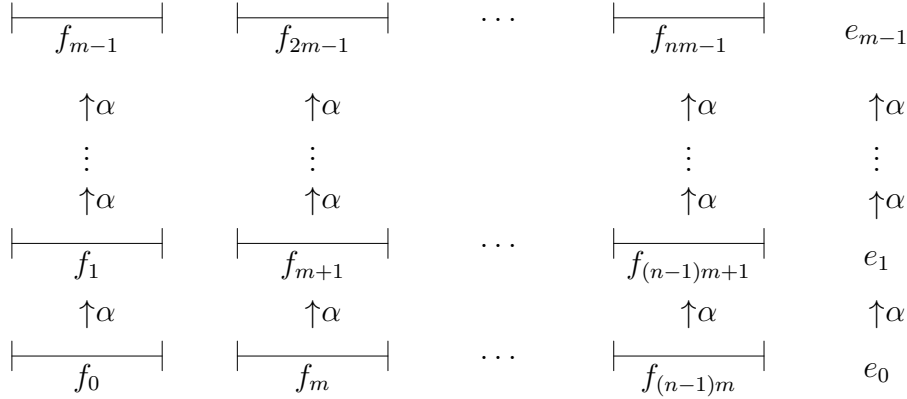

\begin{remark}\label{rem-collapse}From a stack $f_0,f_1,\dots,f_{nm-1}$ of height $nm$ we can construct a stack\linebreak $e_0,e_1,\dots,e_{m-1}$ of height $m$ having the same sum by taking $e_j=\sum_{i=0}^{n-1}f_{im+j},$  and if $\alpha(f_{nm-1})=f_0$ then $\alpha(e_{m-1})=e_0.$  See Figure~\ref{fig-collapse}.
\end{remark}

For a compact space $X$ we write $C(X)$ for the $C^*$-algebra of 
complex-valued continuous functions on $X$ under the supremum norm.  
The unital $C^*$-subalgebra generated by $a\in A$ is denoted $C^*(a).$
We need a few basic approximation results for $C^*$-algebras.  The 
proofs use standard techniques in spectral theory.

\begin{lemma}If $e$ and $f$ are projections in $A$ with $\|ef\|<1/4$ then 
there is a projection $g$ such that $e\perp g$ and $\|f-g\|\leq4\|ef\|.$
\end{lemma}
\begin{proof}
With $e,f$ as in the hypothesis, $(1-e)f(1-e)$ is self-adjoint and close to $f,$ 
\[
\|(1-e)f(1-e)-f\|=\|ef(e-1)-fe\|\leq\|ef\|\|e-1\|+\|fe\|=2\|ef\|.
\]
Since $\spec(f)\subset\{0,1\},$ we have 
$\spec((1-e)f(1-e))\subset B_{2\|ef\|}(0)\bigcup B_{2\|ef\|}(1),$ where 
$B_\delta(x)$ denotes the ball of radius $\delta$ centered at $x.$  The 
function $\varphi$ on $\spec((1-e)f(1-e))$ which is $0$ on $B_{2\|ef\|}(0)$ 
and $1$ on $B_{2\|ef\|}(1)$ is a projection in $C(\spec((1-e)f(1-e)))$ with 
$\sup|\varphi(x)-x|<2\|ef\|,$ so the continuous functional calculus grants 
us a projection $g\in C^*((1-e)f(1-e))$ with $\|g-(1-e)f(1-e)\|<2\|ef\|.$  
Now $eg=0$ (since $e(1-e)=0$) and  
$\|f-g\|\leq\|f-(1-e)f(1-e)\|+\|(1-e)f(1-e)-g\|<4\|ef\|$ as desired.
\end{proof}

\begin{lemma}\label{lem-bumppartial}If $e,f$ are projections in $A$ with 
$\|e-f\|<1/2$ then there is a unitary $u\in A$ with $\|1-u\|\leq 4\|e-f\|$ 
and $u^*eu=f.$
\end{lemma}
\begin{proof}Set $z=2ef-e-f+1.$  One readily checks that $ez=ef=zf,$ and 
$z$ is invertible since $\|z-1\|=\|e(f-e)+(e-f)f\|\leq2\|e-f\|<1.$  
Consider the polar decomposition $z=ua$ where $u$ is a unitary and 
$a=(z^*z)^{1/2}.$  We have $fz^*z=(zf)^*z=(ez)^*z=z^*ez=z^*zf,$ i.e., 
$z^*z$ commutes with $f,$ and since $a\in C^*(z^*z)$ we see that $a$ also 
commutes with $f.$  Thus $u^*eu=az^{-1}eza^{-1}=afa^{-1}=f.$  Finally, 
\[
\|1-u\|\leq\|1-z\|+\|z-u\|=\|1-z\|+\|1-a\|\leq\|1-z\|+\|1-z\|\leq4\|e-f\|.
\]
\end{proof}

\begin{lemma}\label{lem-ortho}If $\{e_1,e_2,\dots,e_n\}$ and 
$\{f_1,f_2,\dots,f_n\}$ are families of orthogonal projections in $A$ with 
$\|e_i-f_i\|<\varepsilon<\frac{1}{2n}$ for each $i$ then there is a unitary 
$u\in A$ such that $\|1-u\|<8n\varepsilon$ and $u^*e_iu=f_i$ for each $i.$
\end{lemma}
\begin{proof}Let $e_0=1-\sum_{i=1}^ne_i$ and $f_0=1-\sum_{i=1}^nf_i.$  By 
the triangle inequality, $\|e_0-f_0\|\leq n\varepsilon.$  By 
Lemma~\ref{lem-bumppartial} we can for each $i$ choose a unitary $u_i$ 
with $\|1-u_i\|<4\|e_i-f_i\|$ such that $u_i^*e_iu_i=f_i.$  Set 
$u=\sum_{i=0}^ne_iu_i.$  Then $u$ is a unitary with $u^*e_iu=f_i$ for each 
$i$ and
\[
\|1-u\|\leq\sum_{i=0}^n\|e_i-u_ie_i\|\leq\sum_{i=0}^n\|1-u_i\|\leq8n\varepsilon.
\]
\end{proof}

\begin{lemma}\label{lem-partial}If $e$ and $f$ are projections in $A$ and 
$p\in A$ is a partial isometry with $\|p^*p-e\|<\varepsilon<1/2$ and 
$\|pp^*-f\|<\varepsilon$ then there is a partial isometry $q\in A$ with 
$q^*q=e$ and $qq^*=f$ and $\|p-q\|<8\varepsilon.$
\end{lemma}

\noindent 
When $e,f\in A$ are equivalent projections we write $e\sim f;$ i.e.,
$e\sim f$ if there exists $p\in A$ such that $p^*p=e$ and $pp^*=f.$
We state two facts about equivalence classes and $K_0(A)$ and refer
the reader to \cite{Bl} for proofs.

\begin{lemma}\label{lem-allproj}If $A$ has cancellation then for any 
$g\in K_0(A)_+$ with $g<[1]$ there is a projection $e\in A$ with $[e]=g.$
\end{lemma}

\begin{definition}When $K_0(A)$ is torsion-free, the rank of $K_0(A)$ is 
the dimension of the rational vector space $K_0(A)\otimes\rats.$
\end{definition}

\begin{lemma}\label{lem-piproj}Suppose $A$ has cancellation.  Write $[\cdot]$ 
for equivalence class in $K_0(A).$  If $e$ and $f$ are projections in $A$ such 
that $[e]\leq[f]$ then there is a partial isometry $p\in A$ with $pp^*\leq f$ 
and $p^*p=e,$ and $pp^*=f$ if $[e]=[f].$
\end{lemma}

Many proofs that a certain automorphism has the Rohlin property use an 
intermediate concept similar to the following, which comes from \cite{Ki1}. 
\begin{definition}\label{def-arp}An automorphism $\alpha$ of a unital $AF$ 
algebra $A$ has the approximate Rohlin property if for any $m,n\in\nats,$ 
any $\varepsilon>0,$ and any finite-dimensional $C^*$-subalgebra $B$ of $A$ 
there is an $\varepsilon$-approximate cyclic stack $e_0,e_1,\dots,e_{m-1}$ 
in $A\bigcap B'$ such that 
\begin{enumerate}
\item $[e_0]=[e_1]=\cdots=[e_{m-1}]$ and
\item $\displaystyle[e_0]\geq n\left[1-\sum_{i=0}^{m-1}e_i\right],$
\end{enumerate}
where $[\cdot]$ denotes the equivalence class in $K_0(A\bigcap B').$
\end{definition}

\begin{remark}\label{rem-arpnoncyclic}Setting $n=1$ in 
Definition~\ref{def-arp} changes nothing since for $n=\ell$ we can start with 
a $\frac{\varepsilon}{\ell}$-approximate stack of height $\ell m$ satisfying 
the definition for $n=1$ and then the construction in 
Remark~\ref{rem-collapse} yields an $\varepsilon$-approximate stack of height 
$m$ satisfying the definition for $n=\ell.$  We can further simplify the 
definition by dropping the requirement that the approximate stack be cyclic, 
as we shall see later.
\end{remark}

\section{General results on the Rohlin property}\label{sec-gen}
\begin{theorem}\label{arp-rp}
An automorphism of unital AF algebra having the approximate Rohlin property 
also has the Rohlin property.
\end{theorem}

\begin{remark} The approximate Rohlin property was introduced in \cite{BEK}, 
in the context of the one-sided shift on the CAR algebra 
$\displaystyle\bigotimes^\infty M_2.$  Definition~\ref{def-arp} generalizes 
the idea to arbitrary AF algebras.  In \cite{Ki1} it is shown (Theorem~4.1) 
that the approximate Rohlin property implies the Rohlin property for unital 
simple AF algebras when $K_0$ is of finite rank and has no infinitesimal 
elements.  In this setting the Rohlin property also implies the approximate 
Rohlin property, but we do not know whether the converse of
Theorem~\ref{arp-rp} holds in general.
\end{remark}

\noindent{\em Proof of Theorem~\ref{arp-rp}. }
Given $m,\varepsilon$ and $F$ we shall use the approximate Rohlin property to 
construct a partition of unity by projections 
$e_{0,0},\dots,e_{0,m-1},e_{1,0},\dots,e_{m}$ satisfying 
Definition~\ref{def-rp}.  The technique is similar to that in the proof of 
Lemma~4.4 in \cite{Ki1} in which the additional assumptions (that $K_0(A)$ be 
of finite rank and have no infinitesimals) are used.  There is also an error 
or omission in the proof there which we have corrected.

We shall assume that the elements of $F$ are of norm $1$ as there is no loss 
of generality in doing so.  The construction depends on parameters 
$n,\ell\in\nats$ (which determine heights of intermediate stacks) and 
$\epsilon>0.$  We require that $n\equiv1\mod m+1,$ $\ell>4,$ and 
$\epsilon<(Nm)^{-3},$ where $N=n\ell.$ 

Let $B$ be a finite dimensional $C^*$-subalgebra of $A$ such that 
$\dist(\alpha^{-k}(x),B)<(Nm)^{-3}$ for all $x\in F, k=0,1,\dots,Nm-1.$  Let 
$e_0,\dots,e_{Nm-1}$ be an $\epsilon$-approximate cyclic stack in 
$A\bigcap B'$ such that
\[
[e]\leq[e_0]=[e_1]=\dots=[e_{Nm-1}],
\]
where $[\cdot]$ denotes equivalence class in $K_0(A\bigcap B')$ and
\[
e=1-\sum_{i=0}^{Nm-1}e_i.
\]
By Lemma~\ref{lem-piproj} there is a partial 
isometry $p_0$ in $A\bigcap B'$ such that $p_0^*p_0=e$ and $p_0p_0^*<e_0.$  
We have 
\[
\|\alpha(p_0)^*\alpha(p_0)-e\|=\|\alpha(e)-e\|<Nm\epsilon
\]
and since
\[
\alpha(p_0)\alpha(p_0)^*<\alpha(e_0),
\]
$\alpha(p_0)\alpha(p_0)^*$ is within $\varepsilon$ of some subprojection of 
$e_1.$  By Lemma~\ref{lem-partial} there is a partial isometry $p_1\in A$ 
such that $p_1^*p_1=e$ and $p_1p_1^*<e_1$ and 
$\|\alpha(p_0)-p_1\|\ll(Nm)^{-3},$ provided $\epsilon$ is chosen sufficiently 
small.  Continuing inductively, assuming $\epsilon$ is small enough, we can 
find partial isometries $p_k,$ $k=0,1,\dots,Nm-1,$ 
such that $p_k^*p_k=e$ and $p_kp_k^*<e_k$ for each $k$ and
\[
\|\alpha(p_k)-p_{k+1}\|<(Nm)^{-3}\quad\text{for } k=0,1,\dots Nm-2.
\]
If $x\in F$ and $k\in\{0,1,\dots,Nm-1\}$ then
\begin{align*}
\|[p_k,x]\|
&\leq\|p_kx-\alpha^k(p_0)x\|+\|[\alpha^k(p_0),x]\|+\|x\alpha^k(p_0)-xp_k\|\\
&\leq2\|p_k-\alpha^k(p_0)\|+\|[p_0,\alpha^{-k}(x)]\|\\
&\leq2k(Nm)^{-3}+2\dist(\alpha^{-k}(x),B)\\
&<2(Nm)^{-2}.
\end{align*}

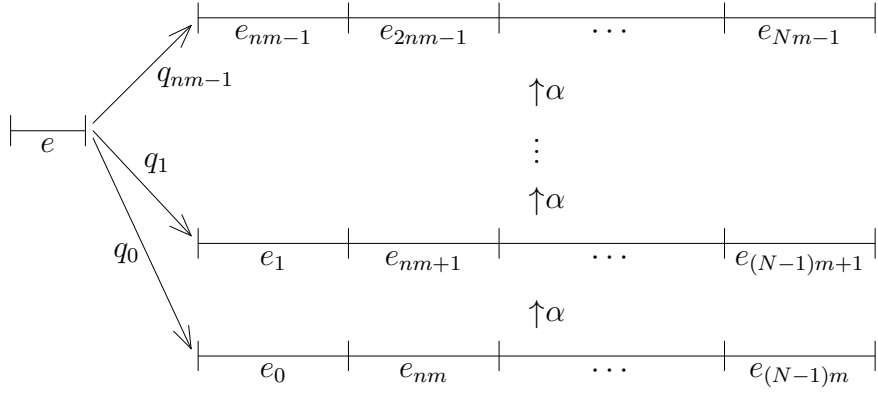
\begin{figure}
\setlength{\unitlength}{1mm} {
\renewcommand{\dashlinestretch}{2}
\renewcommand{\drawlinestretch}{1000}
\begin{center}
\begin{picture}(120,55)(-28,-6)

\drawline(-25,30)(-15,30)
\drawline(-25,28)(-25,32)
\drawline(-15,28)(-15,32)
\put(-20,29){\makebox(0,0)[ct]{$e$}}

\drawline(-14,29)(-1,1)
\drawline(-1,1)(-1.3,4.1)
\drawline(-1,1)(-3.3,3.3)
\put(-9.5,15){\makebox(0,0)[ct]{$q_0$}}

\drawline(-14,30)(-1,16)
\drawline(-1,16)(-3.8,17.3)
\drawline(-1,16)(-2.2,18.8)
\put(-5.5,27){\makebox(0,0)[ct]{$q_1$}}

\drawline(-14,31)(-1,44)
\drawline(-1,44)(-2.3,41.1)
\drawline(-1,44)(-3.9,42.6)
\put(-5.5,37){\makebox(0,0)[l]{$q_{nm-1}$}}

\drawline(0,0)(90,0)
\drawline(0,-2)(0,2)
\drawline(20,-2)(20,2)
\drawline(40,-2)(40,2)
\drawline(70,-2)(70,2) 
\drawline(90,-2)(90,2)

\put(10,-1){\makebox(0,0)[ct]{$e_0$}}
\put(30,-1){\makebox(0,0)[ct]{$e_{nm}$}}
\put(55,-1){\makebox(0,0)[ct]{$\cdots$}}
\put(80,-1){\makebox(0,0)[ct]{$e_{(N-1)m}$}}

\put(45,7.5){\makebox(0,0)[ct]{$\uparrow$}}
\put(47.5,6.5){\makebox(0,0)[ct]{$\alpha$}}

\drawline(0,15)(90,15)
\drawline(0,13)(0,17)
\drawline(20,13)(20,17)
\drawline(40,13)(40,17)
\drawline(70,13)(70,17) 
\drawline(90,13)(90,17)

\put(10,14){\makebox(0,0)[ct]{$e_1$}}
\put(30,14){\makebox(0,0)[ct]{$e_{nm+1}$}}
\put(55,14){\makebox(0,0)[ct]{$\cdots$}}
\put(80,14){\makebox(0,0)[ct]{$e_{(N-1)m+1}$}}

\put(45,22.5){\makebox(0,0)[ct]{$\uparrow$}}
\put(47.5,21.5){\makebox(0,0)[ct]{$\alpha$}}

\put(45,31){\makebox(0,0)[ct]{$\vdots$}}

\put(45,37){\makebox(0,0)[ct]{$\uparrow$}}
\put(47.5,36){\makebox(0,0)[ct]{$\alpha$}}

\drawline(0,45)(90,45)
\drawline(0,43)(0,47)
\drawline(20,43)(20,47)
\drawline(40,43)(40,47) 
\drawline(70,43)(70,47)
\drawline(90,43)(90,47)

\put(10,44){\makebox(0,0)[ct]{$e_{nm-1}$}}
\put(55,44){\makebox(0,0)[ct]{$\cdots$}}
\put(30,44){\makebox(0,0)[ct]{$e_{2nm-1}$}}
\put(80,44){\makebox(0,0)[ct]{$e_{Nm-1}$}}

\end{picture}
\end{center}
}
\caption{The situation in Theorem~\ref{arp-rp}.}
\label{fig-arprp}
\end{figure}

\noindent For $k=0,1,\dots, mn-1$ the partial isometry 
\[
q_k=\frac{1}{\sqrt{\ell}}\sum_{j=0}^{\ell-1}p_{jmn+k}
\]
has domain $e$ and range $<\sum_{j=0}^{\ell-1}e_{jmn+k}.$  By the triangle 
inequality, 
\[
\|\alpha(q_k)-q_{k+1}\|<\sqrt{\ell}(Nm)^{-2}< 3/\sqrt{\ell}
\quad k=0,1,\dots,m-2,
\]
\[
\|\alpha(q_{mn-1})-q_{0}\|
<\frac{1}{\sqrt{\ell}}
\left[(\ell-1)(Nm)^{-2}+2\right]<3/\sqrt{\ell},
\text{ and}
\]
\[
\|[q_k,x]\|<2\sqrt{\ell}(Nm)^{-2},\quad x\in F,\ k=0,1,\dots,m-1.
\]
Define $e_{0,k}$ for $k=0,1,\dots,m-1$ by setting 
\[
e_{0,k}=\sum_{j=0}^{N-1}e_{jm+k}-\sum_{j=0}^{n-1}q_{jm+k}q_{jm+k}^*.
\]
Another application of the triangle inequality, together with the bounds 
above, yields
\[
\|\alpha(e_{0,k})-e_{0,k+1}\|< N(Nm)^{-3} + 
6n/\sqrt{\ell}<7n/\sqrt{\ell}\quad k=0,1,\dots,m-1,
\]
where $e_{0,m}=e_{0,0},$ and 
\[
\|[e_{0,k},x]\|<2N(Nm)^{-3}+4n\sqrt{\ell}(Nm)^{-2}<5/N,
\quad x\in F, k=0,1,\dots,m-1,
\]
where we have used the fact that $\|[e_i,x]\|<2(Nm)^{-3}.$ 

We claim that the $C^*$-subalgebra $D$ generated by 
$\{q_0,q_1,\dots,q_{mn-1}\}$ is isomorphic to $M_{mn+1}.$  Indeed, 
$D$ has $mn+1$ orthogonal projections, 
$e,\,q_0q_0^*,\,q_1q_1^*,\,\dots,\,q_{mn-1}q_{mn-1}^*,$ and for each 
$i,$ $q_i$ is a partial isometry from $e$ to $q_iq_i^*.$  Thus for each 
pair $i,j,$ $q_iq_j^*$ is a partial isometry from $q_jq_j^*$ to $q_iq_i^*,$ 
and we see that the generators for $D$ are part of a full set of matrix units 
(which they generate) for an $mn+1$-dimensional $C^*$-algebra.  Let $u$ be 
the unitary in $D$ defined by 
\[
(\Ad u)(q_k)=q_{k+1}
\]
for $k=0,1,\dots,mn-1,$ where $q_{mn}=q_0,$ and 
\[
(\Ad u)(e)=e.
\]
From the triangle inequality we obtain the generous estimate
\[
\|\alpha|_D-(\Ad u)|_D\|<3(mn+1)^2/\sqrt{\ell}.
\]
The eigenvalues of $u$ are 
$e^{\frac{2\pi ij}{mn}},$ $j=0,1,\dots mn,$ and
there is a unitary $v\in D$ with eigenvalues $e^{\frac{2\pi ij}{mn+1}},$ 
$j=0,1,\dots mn,$ such that $\|u-v\|<\frac{2\pi}{nm}.$  To see this, 
note that unitaries in $M_{mn+1}$ can be diagonalized by unitaries in 
$M_{mn+1}$ and there is a path of diagonal unitaries from 
$\diag(1,\frac{2\pi i}{mn},\dots,\frac{2\pi i(mn-1)}{mn},1)$ to 
$\diag(1,\frac{2\pi i}{mn+1},\dots,\frac{2\pi i(mn-1)}{mn+1},\frac{2\pi i(mn)}{mn+1})$ 
having length $\frac{2\pi}{nm+1}.$  There is a projection $r$ in $D$ such that 
$\sum_{j=0}^{mn}(\Ad v)^j(r)$ equals the identity in $D$ ($v$ acts as a permutation 
of least period $mn+1$ on some orthonormal basis for $\cplx^{mn+1};$ 
we can take $r$ to be the projection onto one of these basis vectors) and we set 
\[
e_{1,k}=\sum_{j=0}^{\frac{nm+1}{m+1}-1}(\Ad v)^{j(m+1)}(r), \quad k=0,1,\dots,m,
\]
which makes sense because $n\equiv1\mod m+1.$ For $k=0,1,\dots,m$ we have 
$(\Ad v)(e_{1,k})=e_{1,k+1}$ and 
\begin{align*}
\|\alpha(e_{1,k})-e_{1,k+1}\|&\leq\|\alpha(e_{1,k})-(\Ad u)(e_{1,k})\|+\|(\Ad u)(e_{1,k})-(\Ad v)(e_{1,k})\|\\
&\leq 3(mn+1)^2/\sqrt{\ell}+\frac{4\pi}{nm},
\end{align*}
where $e_{1,m+1}=e_{1,0}.$  For $x\in F$ we bound $\|[e_{1,k},x]\|$ by 
$4(nm+1)^2\sqrt{\ell}(Nm)^{-2}<4/\sqrt{\ell},$ using the triangle inequality 
with our previous estimates.

Now $\{e_{i,j}:\ i=0,1,\ j=0,1,\dots,m+i-1\}$ is a partition of unity by 
projections.  Our estimates show that we achieve the desired tightness in 
our approximate Rohlin stacks if we take $n>\frac{8\pi}{m\epsilon}$ and 
$\ell>\frac{36(mn+1)^4}{\epsilon^2}.$
\qed

Inspection of the proof shows that one can deduce the Rohlin property from 
weaker hypotheses, by simply replacing the hypotheses involving $K_0$ in the 
approximate Rohlin property with a new requirement that a certain partial 
isometry exists.  This implies an inequality in $K_0$ but is not implied by 
it for algebras without nice properties like cancellation.
\begin{theorem}\label{nokrohlin}Let $\alpha$ be an automorphism of a unital 
$C^*$-algebra $A$ with the property that for any finite subset $F$ of $A$ and 
for any $m\in\nats,$ $\epsilon>0,$ there are projections 
$e_0,e_1,\dots,e_{m-1}$ in $A$ and a partial isometry $p$ in $A$ such that 
\begin{enumerate}
\item $\|e_ie_j\|<\epsilon$ for $i,j=0,1,\dots,m-1,$ $i\neq j,$
\item $\|\alpha(e_i)-e_{i+1}\|<\epsilon$ for $i=0,1,\dots,m-1,$ where $e_m=e_0,$
\item $\displaystyle\left|p^*p-\left(1-\sum_{i=0}^{m-1}e_i\right)\right|<\epsilon$ and $\|pp^*-e_0\|<\epsilon,$
\item $\|[x,y]\|<\epsilon$ for any $x\in F, y\in\{e_0,p\}.$
\end{enumerate}
Then $\alpha$ has the Rohlin property.
\end{theorem}
\begin{proof}
We have assumed what we had to deduce from $K$-theoretic data in the proof of 
Theorem~\ref{arp-rp}, and we repeat that proof with minor changes.   Use 
$\displaystyle\bigcup_{i=0}^{-Nm+1}\alpha^{-i}(F)$ in place of $F.$ Perturb 
the $e_i$ so that they are orthogonal (Lemma~\ref{lem-ortho}).  Let $p_0$ be 
a partial isometry close to $p$ for which $p_0^*p_0=1-\sum_{i=0}^{Nm-1}e_i$ 
and $p_0p_0^*=e_0$ (Lemma~\ref{lem-piproj}).  The rest of the proof is the 
same.
\end{proof}

We end the section with a condition which is sometimes easier to verify, to 
be used in the proof of Theorem~\ref{sft-rp}.  The hypothesis that the top 
of the stack be taken near the bottom is removed at the expense of a requiring 
a partial isometry between the first two levels of the stack.  This implies 
equality of the $[e_i]$ always, but is not implied by it in algebras without 
cancellation.
\begin{theorem}\label{prohlin}Let $\alpha$ be an automorphism of a unital 
$C^*$-algebra $A$ with the property that for any finite subset $F$ of $A$ 
and for any $m\in\nats,$ $\varepsilon>0,$ there are an 
$\varepsilon$-approximate stack $e_0,e_1,\dots,e_{m-1}$ for $\alpha$ and 
partial isometries $p,q$ in $A$ such that 
\begin{enumerate}
\item $p^*p=1-\sum_{i=0}^{m-1}e_i$ and $pp^*<e_0,$
\item $q^*q=e_0$ and $qq^*=e_1,$
\item $\|[x,y]\|<\varepsilon$ for any $x\in F, y\in\{p,q\}.$
\end{enumerate}
Then $\alpha$ has the Rohlin property.
\end{theorem}
\begin{remark}We could allow $e_i$ to be almost orthogonal instead of 
orthogonal and relax the equalities involving $p$ and $q$ to approximations 
as in Theorem~\ref{nokrohlin} and the conclusion still holds.
\end{remark}
\begin{proof}
We show that $\alpha$ satisfies the hypothesis of Theorem~\ref{nokrohlin}.  
Our proof uses many of the same ideas as in  Lemmas~2.1 and 4.3 in \cite{Ki1}.  
The construction depends on parameters 
$\ell\in\nats,\epsilon>0$ to be specified later.  Let $F$ be a finite set, 
$m\in\nats,$ and $\epsilon>0.$ Let $e_0,e_1,\dots,e_{(\ell-1)(m+2)m-1}$ be an 
$\varepsilon$-approximate stack for $\alpha$ and let $p,q$ be partial 
isometries in $A$ such that
\begin{itemize}
\item $\displaystyle p^*p=1-\sum_{i=0}^{(\ell-1)(m+2)m-1}e_i$ and $pp^*<e_0$,
\item $q^*q=e_0$ and $qq^*=e_1$,
\item $\|[x,y]\|<\varepsilon$ for all 
$\displaystyle x\in\bigcup_{i=0}^{(\ell-1)(m+2)}\alpha^{-i}(F),$ $y\in\{p,q\}.$
\end{itemize}
By Lemma~\ref{lem-ortho} there is a unitary $u\in A$ such that 
\[
\|1-u\|<8(\ell-1)(m+2)m\varepsilon
\]
and
\[
u\alpha(e_i)u^*=e_{i+1},\quad i\in\{0,1,\dots,(\ell-1)(m+2)m-2\}.
\]
Write $\beta$ for $(\Ad u)\circ\alpha.$  Define a system of matrix units 
$\{q_{i,j}:0\leq i,j \leq (\ell-1)(m+2)m-1\}$ by setting 
\[
q_{i,i}=e_i
\]
and for $i<j$
\[
q_{i,j}=\beta^i(q)\beta^{i+1}(q)\cdots\beta^{j-1}(q),
\quad q_{j,i}=q_{i,j}^*.
\]
One readily checks that 
\[
q_{i,j}q_{i',j'}=\delta_{j,i'}q_{i,j'},
\]
and
\[
\beta(q_{i,j})=q_{i+1,j+1},\quad 0\leq i,j\leq (\ell-1)(m+2)m-2\}.
\]
For $j=0,1,\dots,m-1$ put 
\begin{align*}
f_j=&\sum_{i=0}^{\ell-2}\bigg[\frac{i+1}{\ell}q_{im+j,im+j}+
\frac{\ell-1-i}{\ell}q_{((m+1)(\ell-1)+i)m+j,((m+1)(\ell-1)+i)m+j})\\
	&\qquad+\frac{\sqrt{(i+1)(\ell-1-i)}}{\ell}q_{im+j,((m+1)(\ell-1)+i)m+j}\\
	&\qquad+\frac{\sqrt{(i+1)(\ell-1-i)}}{\ell}q_{((m+1)(\ell-1)+i)m+j,im+j}\bigg]\\
	&+\sum_{i=\ell-1}^{(m+1)(\ell-1)-1}q_{im+j,im+j}
\end{align*}
Then $f_0,f_1,\dots,f_{m-1}$ is a stack of height $m$ for $\beta$.  
We have {\small
\begin{align*}
\beta(f_{m-1})-f_0=
&\frac{1}{\ell}\sum_{i=0}^{\ell-1}
\bigg[(e_{\scriptscriptstyle((m+1)(\ell-1)+i)m}-e_{\scriptscriptstyle im})\\
&\quad+\left(\sqrt{i(\ell-i)}-\sqrt{(i+1)(\ell-1-i)}\right)
(q_{\scriptscriptstyle im,((m+1)(\ell-1)+i)m}+
q_{\scriptscriptstyle((m+1)(\ell-1)+i)m,im})\bigg],
\end{align*}
}
and since, for $i=0,1,\dots,\ell-1,$
\begin{align*}
\left|\sqrt{i(\ell-i)}-\sqrt{(i+1)(\ell-1-i)}\right|
&\leq\sqrt{|i(\ell-i)-(i+1)(\ell-1-i)|}\\
&=\sqrt{|1+2i-\ell|}\\
&<\sqrt{\ell}
\end{align*}
we see from orthogonality that 
\[
\|\beta(f_{m-1})-f_0\|<\frac{1}{\ell}+\frac{1}{\sqrt{\ell}}.
\]
This shows that $f_0,f_1,\dots,f_{m-1}$ is a 
$\frac{2}{\sqrt{\ell}}$-approximate cyclic stack for $\beta,$ and thus
by choosing $\varepsilon$ small enough and $\ell$ large enough we can ensure 
that $f_0,f_1,\dots,f_{m-1}$ is an $\epsilon$-approximate cyclic stack for 
$\alpha.$  We have constructed an approximate cyclic stack from an 
approximate noncyclic stack, and the sum of this new stack is smaller than 
the sum of the old one by $m(\ell-1)$ projections equivalent to $e_0,$ 
specifically, 
\begin{align*}
1-\sum_{j=0}^{m-1}f_j=p^*p+\sum_{j=0}^{m-1}\sum_{i=0}^{\ell-2}\bigg\{
&\frac{\ell-1-i}{\ell}q_{im+j,im+j}\\
&+\frac{i+1}{\ell}q_{((m+1)(\ell-1)+i)m+j,((m+1)(\ell-1)+i)m+j}\\
&+\frac{\sqrt{(i+1)(\ell-1-i)}}{\ell}q_{im+j,((m+1)(\ell-1)+i)m+j}\\
&+\frac{\sqrt{(i+1)(\ell-1-i)}}{\ell}q_{((m+1)(\ell-1)+i)m+j,im+j}\bigg\}.
\end{align*}
The base of our new stack, $f_0$, has a subprojection equivalent to 
$\displaystyle 1-\sum_{j=0}^{m-1}f_j,$ indeed
\begin{align*}
r=\sum_{j=0}^{m-1}\sum_{i=0}^{\ell-2}&\left(\sqrt{\frac{\ell-1-i}{\ell}}q_{(\ell+im+j-1)m,im+j}
+\sqrt{\frac{i+1}{\ell}}q_{(\ell+im+j-1)m,((m+1)(\ell-1)+i)m+j}\right)\\
&+\left(\sqrt{\frac{1}{\ell}}e_0+\sqrt{\frac{\ell-1}{\ell}}q_{(m+1)(\ell-1)m,0}\right)p
\end{align*}
is a partial isometry satisfying
\[
rr^*=1-\sum_{j=0}^{m-1}f_j\quad\text{and}\quad r^*r<f_0.
\]
The $q_{i,j}$ almost commute with $F$ if $\epsilon$ is chosen small relative to 
$1/\ell,$ and since $f_0,\dots,f_{m-1}$ and $r$ are linear combinations of the 
$q_{i,j}$ these too almost commute with $F.$  By Theorem~\ref{prohlin}, 
$\alpha$ has the Rohlin property.
\end{proof}

\begin{remark}In the proof of Theorem~\ref{prohlin} we constructed an 
approximate cyclic stack $f_0,f_1,\dots,f_{m-1}$ from an approximate 
(noncyclic) stack $e_0,e_1,\dots,e_{Nm-1}$ using only the partial isometry 
$q$ from $e_0$ to $e_1.$  To construct the partial isometry $r$ from 
$1-\sum f_j$ to a subprojection of $f_0$ we needed the partial isometry 
$p$ from $1-\sum e_i$ to a subprojection of $e_0.$  The $f_j$ and $r$ lie in 
$C^*(p,q,\alpha(q),\dots,\alpha^{Nm-1}(q))$ and the proof shows that these can 
be made to almost commute with any finite set provided the $e_i$ and $q$ and 
$p$ can be made to almost commute with any finite set.  This justifies the 
last line of Remark~\ref{rem-arpnoncyclic}, since if we can make the $f_j$ 
and $r$ almost commute with a system of matrix units for $B$ then a small 
inner perturbation takes the $f_j$ and $r$ into the commutant of $B.$
\end{remark}

\section{Shifts of finite type and the bilateral algebras}\label{sec-def} 

We refer the reader to the excellent book of Lind and Marcus \cite{LM} 
for background on shifts of finite type.  Let $T$ be a nonnegative 
integral $r\times r$ 
matrix.  Let $G_T$ be the directed graph with $r$ vertices labeled 
$1,2,\dots,r$ and having $T(i,j)$ edges from vertex $i$ to vertex $j.$  Write 
$\Edge_T$ for the edge set of $G_T,$ and for $e\in\Edge_T$ let $\iver(e)$ and 
$\tver(e)$ be the initial and terminal vertices of $e,$ respectively.

\begin{definition}A path on $G_T$ is a finite or infinite sequence $x$ in 
$\Edge_T$ such that $\iver(x_{i+1})=\tver(x_i)$ for all relevant indices $i.$
\end{definition}

\noindent{\bf Notation. }For a finite path $x=x_0x_1\dots x_{\ell-1}$ we let 
$\iver(x)=\iver(x_0)$ and $\tver(x)=\tver(x_{\ell-1})$ and we write $|x|=\ell$ 
for the length of $x.$  By convention, each vertex $i$ is regarded as a path 
of length $0$ having $i$ as its initial and terminal vertices.  If $a$ and 
$b$ are integers with $a\leq b$ and $x$ is a path then we put 
\[
x_{[a,b]}=x_ax_{a+1}\dots x_b
\]
when the right-hand side makes sense.

We use the following fact often.
\begin{lemma}\label{lem-pathcount}For $i,j\in\{1,2,\dots,r\},$ the number of 
paths of length $\ell$ in $G_T$ starting at $i$ and ending at $j$ is 
$T^{\ell}(i,j).$
\end{lemma}

\begin{definition}The underlying space $X_T$ of the shift of finite type 
determined by $T$ is the set of all bi-infinite paths on $G_T.$
\end{definition}

\begin{definition}
Basic cylinders are sets of the form
\[
\cyl(x,k)=\{y\in X_T:y_{i+k}=x_i\text{ for }0\leq i<|x|\},
\]
where $x=x_0x_1\dots x_{|x|-1}$ is a finite path in $G_T$ and $k\in\ints.$
\end{definition}

\noindent Basic cylinders form a basis of clopen sets for a topology on $X_T$ 
in which $X_T$ is a Cantor set.  Though we shall not specifically use a metric 
it is helpful to think of two points being close if their sequences agree on a 
large interval of indices about $0.$

The paradigm of symbolic dynamics is that all the complexity of a system is 
encoded in sequence space. 
\begin{definition}
The shift homeomorphism $\sigma_T:X_T\to X_T$ is defined by
\[
\left(\sigma_T\left((x_i)_{i\in\ints}\right)\right)_i=x_{i+1}.
\]
\end{definition}

\begin{example}The full $n$-shift is the shift of finite type associated 
with the $1\times1$ matrix $(n),$ the underlying space of which is 
$\{1,2,\dots,n\}^\ints.$
\end{example}

We shall assume our matrix $T$ is primitive, i.e., some positive power of $T$ 
has strictly positive entries, and that $T\neq(1).$  Then  $X_T$ is nonempty 
and $\sigma_T$ is topologically mixing (see, for example, Proposition~4.5.10 
of \cite{LM}).  By the Perron-Frobenius theorem, $T$ has a unique (up to 
scalar multiplication) strictly positive (right) eigenvector; the 
corresponding eigenvalue $\lambda_T$ is algebraically simple, greater than 
$1,$ and strictly larger in absolute value than any other eigenvalue of $T.$  
Let $v,w^*$ be positive Perron-Frobenius eigenvectors for $T$ and $T^*,$ 
respectively, normalized so that $wv=1.$  

The {\em measure of maximal entropy} on $X_T$ is the Markov measure $\mu_T$ 
defined from transition probabilities 
\[
p(e)=\lambda_T^{-1}v_{\tver(e)}/v_{\iver(e)},\quad e\in\Edge_T
\]
and stationary distribution 
\[
\pi_i=w_iv_i,\quad 1\leq i\leq r
\]
by
\begin{align}\label{eqn-trmu}
\mu_T(\cyl(x,k))
&=\pi_{\iver(x)}p(x_0)p(x_{1})\dots p(x_{|x|-1})\\
&=\lambda_T^{-|x|}w_{\iver(x)}v_{\tver(x)}.
\end{align}
We have the following standard result.
\begin{lemma}\label{lem-trmu}$\mu_T$ is an ergodic measure on $(X_T,\sigma_T)$ 
which gives every cylinder positive measure.
\end{lemma}

We are interested in the symmetric, or bilateral, $C^*$ algebras associated to
shifts of finite type.  The earliest references of which we are aware for the 
these $C^*$-algebras are from 1988,  a construction by Wagoner in \cite{Wa1} 
specifically for shifts of finite type, and Ruelle's more general construction 
in \cite{Ru}.  These algebras are implicit in Krieger's earlier $K$-theoretic 
invariants for shifts of finite type, introduced in \cite{Kr} (1982). 

We write $M_n$ for the $C^*$-algebra of complex $n\times n$ matrices.
For integers $a\leq b$ let $H_T^{[a,b]}$ be the $C^*$-algebra generated by 
matrix units $E_{x,y,a,b},$ where $x$ and $y$ are paths in $G_T$ of length 
$b-a+1$ with $\iver(x)=\iver(y)$ and $\tver(x)=\tver(y).$  Matrix units means
\[
E_{x,y,a,b}E_{x',y',a,b}=\delta_{y,x'}E(x,y',a,b)
\]
and
\[
E_{x,y,a,b}^*=E_{y,x,a,b}.
\]

\begin{lemma}
For each interval $I=[a,b]$ we have 
\[
H_T^I\cong\bigoplus_{i,j}M_{T^{|I|+1}(i,j)},
\]
one summand for each pair of vertices, the size of the $ij$th being the number 
of paths of length $|I|+1$ in $G_T$ beginning at vertex $i$ and ending at 
vertex $j.$  
\end{lemma}
\begin{proof}
For each pair of vertices $i,j,$ enumerate the $T^{|I|+1}(i,j)$ paths of length 
$|I|+1$ in $G_T$ with initial vertex $i$ and terminal vertex $j.$  If $x,y$ are 
the $k$th and $\ell$th such paths from $i$ to $j$ then we identify 
$E_{x,y,a,b}$ with the $T^{|I|+1}(i,j)\times T^{|I|+1}(i,j)$ matrix having a 
$1$ in the $k,\ell$th position and zeros everywhere else, and we embed this in 
$\displaystyle\bigoplus_{i,j}M_{T^{|I|+1}(i,j)}$ by taking the direct sum of 
it with matrices of zeros for the other summands.  It is easy to see that this 
induces an isomorphism.
\end{proof}

Let $I=[a,b]$ and $J=[c,d]$ be intervals with $I\subset J.$   Let 
$\psi_{\scriptscriptstyle T,I,J}:H_T^I\to H_T^J$ be the embedding defined for 
matrix units by
\begin{equation}\label{eqn-psidef}
\psi_{\scriptscriptstyle T,I,J}(E_{x,y,a,b})
=\sum\{E_{z,u,c,d}:z_{[a-c,a+b-c]}=x \text{ and } u_{[a-c,a+b-c]}=y\}.
\end{equation}

\begin{lemma}
The matrix of partial multiplicities for $\psi_{\scriptscriptstyle T,I,J}$ 
is the transpose of $(T^*)^{(a-c)}\otimes T^{(d-b)}.$ 
\end{lemma}
\begin{proof}
Simply count the number of paths $z$ and $u$ which extend $x$ and $y,$ 
respectively, using Lemma~\ref{lem-pathcount}.
\end{proof}

\begin{lemma}If $I\subset J\subset K$ then   
\[
\psi_{\scriptscriptstyle T,I,K}
=\psi_{\scriptscriptstyle T,J,K}\circ\psi_{\scriptscriptstyle T,I,J}.
\]
\end{lemma}
\begin{proof}
This follows from the fact that if $z$ is a path extending $u$ and $u$ is 
a path extending $x$ then $z$ is a path extending $x.$
\end{proof}

\begin{definition}$A_T$ is the closure of 
$\varinjlim(\varphi_{\scriptscriptstyle T,I,J}:H_T^I\to H_T^J).$
\end{definition}

The proof of the following lemma is straightforward and we leave it to the
reader.
\begin{lemma}\label{lem-traceorder}
$A_T$ is simple and has a unique trace.  The trace 
$\tau_{\scriptscriptstyle T}$ is given on each $H_T^I$ by 
\[
\tau_{\scriptscriptstyle T}\left(\sum_{x,y} c_{x,y}E_{x,y,a,b}\right)
=\sum_x c_{x,x}\mu_T(\cyl(x,0)),
\]
where the sums are over all paths of length $|I|+1.$  
The trace determines the order on $K_0(A_T)$ up to infinitesimals.
\end{lemma}

The shift $\sigma_T$ on $X_T$ promotes to an automorphism $\alpha_T$ of $A$ 
defined by
\begin{equation}\label{eqn-alphadef}
\alpha_T(E_{u,v,a,b})=E_{u,v,a-1,b-1}.
\end{equation}
Uniqueness of the trace implies that it is invariant under $\alpha_T,$ i.e.,
\[
\tau_{\scriptscriptstyle T}\circ\alpha_T=\tau_{\scriptscriptstyle T}.
\]

The $C^*$-algebra $C(X_T)$ of continuous complex valued functions on $X_T$ is 
a subalgebra of $A_T$ in a natural way.  Specifically, 
$\varphi_T:C(X_T)\to A_T$ is defined for characteristic functions of basic 
cylinders by 
\[
\varphi_T(\chi_{\cyl(u,k)})=E_{u,u,k,k+|u|-1}
\]
and extended to all of $C(X_T)$ in the usual way.  The next two lemmas follow 
from the definition of $\alpha_T$ in Eq. (\ref{eqn-alphadef}).

\begin{lemma}The $C^*$-embedding $\varphi_T$ takes the induced map 
$f\mapsto f\circ\sigma_T^{-1}$ on $C(X_T)$ to $\alpha_T,$ i.e., 
\[
\varphi_T(f\circ\sigma_T^{-1})=\alpha_T(\varphi_T(f)).
\]
\end{lemma}

\begin{lemma}\label{lem-restrict}
The restriction $\alpha_T|H_T^I$ is an isomorphism from $H_T^I$ to $H_T^{I-1}.$
\end{lemma}

\noindent The following is a very useful fact about these algebras.
\begin{lemma}\label{lem-commute}If $I$ and $J$ are disjoint intervals then 
elements of $H_T^I$ commute with elements of $H_T^J$ in $A_T.$
\end{lemma}
\begin{proof}
We may assume that $I=[a,b]$ and $J=[c,d]$ with $b<c.$  Let $K=[a,d].$  
It follows from the definitions that 
\[
\psi_{\scriptscriptstyle T,I,K}(E_{x,y,a,b})
\psi_{\scriptscriptstyle T,J,K}(E_{x',y',c,d})
=\sum_tE_{xtx',yty',a,d}
\]
where the sum is over all paths $t$ from $\tver(x)$ to $\iver(x')$ of length 
$c-b-1$ and where $xtx'$ means the concatenation of the paths $x,t,x',$ and 
similarly for $yty'$.  Nothing here depends on the order of multiplication and 
$\psi_{\scriptscriptstyle T,I,K}(E_{x',y',c,d})
\psi_{\scriptscriptstyle T,J,K}(E_{x,y,a,b})$ is 
easily seen to be the very same sum.
\end{proof}

\section{Shift equivalence}\label{sec-SE}
Shift equivalence and strong shift equivalence for shifts of finite type were 
introduced by Williams in \cite{Wi}.
\begin{definition}\label{def-se}Nonnegative integral square matrices $U,V$ 
are said to be shift equivalent if there exist nonnegative integral matrices 
$R,S$ and a positive integer $\ell$ such that
\begin{align*}
RS=U^\ell&\qquad SR=V^\ell\\
SU=VS&\qquad UR=RV.
\end{align*}
\end{definition}
Shift equivalence is an equivalence relation.  The integer $\ell$ is called 
the lag of the shift equivalence.
\begin{definition}\label{def-sse}An elementary strong shift equivalence 
between nonnegative integral \linebreak square matrices $U$ and $V$  is a pair of 
nonnegative integral matrices $R,S$ such that $RS=U$ and $SR=V.$  Strong 
shift equivalence is the equivalence relation generated by elementary 
strong shift equivalence.
\end{definition}

Shifts of finite type $(X_U,\sigma_U)$ and $(X_V,\sigma_V)$ are {\em conjugate}
if there is a homeomorphism $\rho:X_U\to X_V$ such that 
$\rho\circ\sigma_U=\sigma_V\circ\rho$ and {\em eventually conjugate} if 
$(X_U,\sigma_U^\ell)$ and $(X_V,\sigma_V^\ell)$ are conjugate for all 
sufficiently large $\ell.$  Williams showed (\cite{Wi}) that 
$U$ and $V$ are strong shift equivalent if and only if 
$(X_U,\sigma_U)$ and $(X_V,\sigma_V)$ are conjugate, and that $U$ and $V$ 
are shift equivalent if and only if $(X_U,\sigma_U)$ and $(X_V,\sigma_V)$ 
are eventually conjugate.

The Shift Equivalence Conjecture is the proposition that shift equivalence 
is the same as strong shift equivalence.  The first counterexample was a 
pair of reducible shifts of finite type found in 1992 by Kim and Roush 
(\cite{KR1}).  The question for the irreducible case remained open until 
1997, when Kim and Roush exhibited in \cite{KR2} a pair of primitive 
$7\times7$ matrices which are shift equivalent but not strong shift 
equivalent.  The invariant used to prove the matrices are not shift equivalent 
was the sign-gyration-compatibility-condition class $sgc_2,$ which turns out 
to be the same as Wagoner's $K_2$-invariant, $\Phi_2.$

There are still many open questions.  For one, it is not known whether the 
conjecture holds for matrices shift equivalent to full shifts.  We do not even 
know if the automorphisms $\alpha_U$ and $\alpha_V$ associated with shift 
equivalent matrices $U$ and $V$ are conjugate, though they are of course 
eventually conjugate.  We do have the following.

\begin{proposition}\label{prop-se}If $U$ and $V$ are shift equivalent 
primitive matrices then there is an isomorphism 
$\Psi:A_U\to A_V$ such that $(\alpha_U)_*
=(\Psi^{-1}\circ\alpha_V\circ\Psi)_*.$
\end{proposition}

\noindent{\em Sketch of Proof. }
This can be proved using the standard bipartite graph construction 
(see \cite{LM}, \S7.2) to obtain a conjugacy between $(X_U,\sigma_U^\ell)$ 
and $(X_V,\sigma_V^\ell)$ for some $\ell.$  A conjugacy produced in this 
way induces a isomorphism of $C^*$-algebras with the desired property.
\qed

We now show that our automorphisms are usually not approximately inner.  One 
consequence of this is that we cannot deduce the Rohlin property for 
$\alpha_T$ from the result of \cite{Ki1}.

\noindent{\bf Notation. }For an $s\times s$ integral matrix $S$ we write 
$\erange_S$ for the eventual range of $S,$
\[
\erange_S=\bigcap_{n\geq0}(\rats^s)^*S^n.
\]
By rank considerations, we have $\erange_S=(\rats^s)^*S^s.$  The eventual rank 
of $S$ is the dimension of the rational vector space $\erange_S\otimes\rats.$

\begin{proposition}\label{prop-ian}The following are equivalent:
\begin{enumerate}
\item $K_0(A_T)$ has no infinitesimals 
\item $T$ has eventual rank $1$ 
\item $T$ is shift equivalent to the $1\times1$ matrix $(n)$ for some $n\geq 2$ 
\item $\alpha_T$ is approximately inner.
\end{enumerate}
\end{proposition}
\begin{proof}
{\em 1 $\Leftrightarrow$ 2. }As ordered groups, 
$K_0(A_T)\cong\varinjlim(\ints^r\otimes\ints^r,T^*\otimes T),$ thus
\[
\rk K_0(A_T)= \rk\erange_{T^*\otimes T} = (\rk\erange_T)^2.
\]
We shall use the fact that 
$\varinjlim(\ints^r\otimes\ints^r,T^*\otimes T)$ is isomorphic to the 
(additive) subgroup of $M_r(\rats)$ consisting of those matrices $U$
in $T^rM_r(\rats)T^r$ such that 
$T^nUT^n\in M_r(\ints)$ for some $n\in\nats.$  The 
positive cone is the set of $U$ such that 
$T^nUT^n\in M_r(\ints_+)$ for some $n\in\nats.$  Let us identify 
$(K_0(A_T),\leq)$ with this ordered group.  By Perron-Frobenius theory, 
$U\in K_0(A_T)_+$ if and only if $wUv>0,$ thus the infinitesimals are 
precisely those $U$ for which $wUv=0.$  It is possible to replace $v$ 
with a multiple so that its entries lie in and span 
$\rats[\lambda_T].$ The same is true of $w.$  It follows that we
can find $r^2-\deg\lambda_T$ linearly independent matrices 
$U_i\in M_r(\ints)$ such that $wU_iv=0$ for each $i.$  The infinitesimal 
subgroup of $K_0(A_T)$ is generated by the matrices 
$T^rU_iT^r,$ and its rank is easily seen to be
$(\rk\erange_T)^2-\deg\lambda_T.$  This last quantity is zero if and
only if $\rk\erange_T=1.$

\noindent{\em 2 $\Leftrightarrow$ 3 }is well known.

\noindent{\em 1 $\Rightarrow$ 4. }If $K_0(A_T)$ has no infinitesimals then 
$\alpha_*$ must be the identity on $K_0(A_T)$, since $\alpha_*(g)$ can only 
differ from $g$ by an infinitesimal.

\noindent{\em 4 $\Rightarrow$ 2. }The existence of infinitesimals implies 
that $T$ has an eigenvalue $\lambda$ with $0<|\lambda|<\lambda_T.$  We can
find in the above representation of $K_0(A_T)$ a matrix $U$ with large
entries such that $UT\approx\lambda U$ and $TU\approx\lambda_TU,$ and 
then $\alpha_*(U)\approx\lambda^{-1}\lambda_T(U)\not\approx U.$ 
\end{proof}

\section{The Rohlin property for $\alpha_T$}\label{sec-rp}

\begin{theorem}\label{sft-rp}If the matrix $T$ is primitive then $\alpha_T$ 
has the Rohlin property.
\end{theorem}
Fix a primitive matrix $T.$  To simplify notation we drop the subscript 
$T$ from the various objects defined in the previous sections.  Before proving 
the theorem we need some preliminary results.  

\begin{lemma}\label{asindlem}For any $x,y\in A$ 
\[
\|[\alpha^n(x),y]\|\to 0 \text{ as }|n|\to\infty.
\]
\end{lemma}
\begin{proof}
Given $x$ and $y$ we can find an interval $I$ and elements $x',y'\in H_T^I$ 
such that $\|x-x'\|$ and $\|y-y'\|$ are as small as desired. If $n>|I|$ then 
$I\cap(I+n)=\emptyset$ and since
$\alpha^n(x')\in H_T^{I+n}$ we see from Lemma~\ref{lem-commute} that 
$\alpha^n(x')$ commutes with $y'.$ Thus
\begin{align*}
\|\alpha^n(x)y-y\alpha^n(x)\|
&\leq\|\alpha^n(x)y-\alpha^n(x')y\|+\|\alpha^n(x')y-\alpha^n(x')y'\|\\
&\qquad+\|\alpha^n(x')y'-y'\alpha^n(x')\|+\|y'\alpha^n(x')-y\alpha^n(x')\|\\
&\qquad+\|y\alpha^n(x')-y\alpha^n(x)\|\\
&\leq\|x-x'\|\|y\|+\|x'\|\|y-y'\|+0+\|y'-y\|\|x'\|+\|y\|\|x'-x\|\\
&\leq 2\|y\|\|x-x'\|+2(\|x\|+\|x-x'\|)\|y-y'\|
\end{align*}
for all $n>|I|.$  By making $\|x-x'\|$ and $\|y-y'\|$ sufficiently small we can 
make this last quantity as small as desired for all sufficiently large $n.$
\end{proof}

\begin{remark}\label{rem-asind}Lemma~\ref{asindlem} provides an immediate 
simplification of our situation.  In proving the Rohlin property for $\alpha$ 
we may ignore the almost commuting requirements, since for any finite 
collection of projections and partial isometries satisfying all the other 
properties their images under $\alpha^n,$ for sufficiently large $n,$ have the 
additional property of almost commuting with the finite set $F.$
\end{remark}

Our next lemma says that $\alpha_*$ has arbitrarily small non-infinitesimal 
fixed points.
\begin{lemma}\label{lem-fp}For every $n>0$ there exists $g\in K_0(A)_+$ such that $\alpha_*(g)=g$ and $ng<[1].$
\end{lemma}
\begin{proof}
For any interval $I=[a,b],$ $K_0(H^I)$ is isomorphic to a subspace of
$M_r(\ints)$ via the map that takes the class of a minimal 
projection $E_{u,u,a,b}\in H^I$ to the
matrix having $(\iver(u),\tver(u))$th entry equal to one and all 
other entries equal to zero.
With this identification the class of the unit in $K_0(H^I)$ is
\[
[1]_I=T^{|I|+1}.
\]
Choose $m$ such that $nI<T^m.$  Let $g\in K_0(A)$ be the element 
represented by the identity matrix $I$ in $K_0(H^{[1,m]}).$  It is clear
that $g$ is positive and $ng<1.$  We claim that 
this $g$ is fixed by $\alpha_*.$  Indeed, $g$ is represented by 
$TI$ in $K_0(H^{[0,m]}),$ while $\alpha_*(g)$ is represented by $I$ in 
$K_0(H^{[0,m-1]})$ and thus by $IT$ in $K_0(H^{[0,m]}).$
\end{proof}

Next we show that the fixed points of $\alpha_*$ are order-dense.
\begin{lemma}\label{divlem}For any $m\in\nats$ there exists $g\in K_0(A)$ 
such that $\alpha_*(g)=g$ and $mg<[1]<(m+1)g.$
\end{lemma}
\begin{proof}By Lemma~\ref{lem-fp} we can find $h\in K_0(A)_+$ fixed by 
$\sigma_*$ such that $m(m+1)h<[1].$  Let $N$ be the greatest positive integer 
for which $Nh<[1].$  While $(N+1)h$ may differ from $[1]$ by an infinitesimal 
we certainly have $(N+2)h> [1].$  Set $g=\lfloor N/m\rfloor h.$ 
We obviously have $m\lfloor N/m\rfloor\leq N.$ Write $N=\ell m+k$ with 
$0\leq k\leq m-1.$  Since $N\geq m(m+1)$ we must have $\ell\geq m+1.$  Now
\[
(m+1)\lfloor N/m\rfloor=(m+1)\ell=\ell m +(m+1)\geq \ell m + k + 2 =N+2,
\]
and thus $g$ has the desired properties.
\end{proof}

\begin{lemma}\label{towlem}For any $m\in\nats$ there exists a clopen subset $C$ of $X$ such that
\begin{enumerate}
\item $C,\sigma(C),\dots,\sigma^{m-1}(C)$ are pairwise disjoint,
\item $1-m[\varphi(\chi_C)]<[\varphi(\chi_C)],$ and
\item $\alpha_*([\varphi(\chi_C)])=[\varphi(\chi_C)].$
\end{enumerate}
\end{lemma}
\begin{proof}
Fix $m$ and let $g$ be as in Lemma~\ref{divlem}.  By Lemma~\ref{lem-traceorder} 
we can find $\varepsilon>0$ such that if $U$ is any clopen set in $X$ with 
$\mu U<\varepsilon$ then $[\varphi(\chi_U)]<[1]-mg.$  Let $x$ be a finite path 
in $G$ long enough so that
\[
\mu\left(\bigcup_{k=1}^{m-1}\cyl(x,-k)\right)<\varepsilon/2.
\]
We can find $N\in\nats$ large enough so that
\[
\mu\left(\bigcup_{k=0}^{Nm-1}\cyl(x,k)\right)>1-\varepsilon/2.
\]
Let $C'$ be the set of points $\displaystyle y\in\bigcup_{k=0}^{Nm-1}\cyl(x,k)$ 
for which the least $k\geq0$ such that $y\in\cyl(x,k)$ satisfies 
$k\equiv m-1\mod m.$  The sets $\sigma^{k}(C'), 0\leq k\leq m-1$ are pairwise 
disjoint and
\[
\left(\bigcup_{k=0}^{Nm-1}\cyl(x,k)\right)
\setminus\left(\bigcup_{k=1}^{m-1}\cyl(x,-k)\right)
\subset\bigcup_{k=0}^{m-1}\sigma^k(C')
\]
whence
\[
\mu\left(X\setminus\left(\bigcup_{k=0}^{m-1}\sigma^k(C')\right)\right)
<\varepsilon.
\]  
This implies that $[\varphi(\chi_{C'})]>g.$  

By Lemma~\ref{lem-allproj}, $[\varphi(\chi_{C'})]-g$ is the class of some 
projection $e\in A_T$ and by Lemma~\ref{lem-piproj}, $e$ is equivalent to 
some projection $f<\varphi(\chi_{C'}).$  Every projection in $A$ is equivalent 
to one of the form $\varphi(\chi_E)$ for some clopen set $E$ (this is just the 
statement that every projection in a finite direct sum of full matrix algebras 
is equivalent to a diagonal matrix), so we can find a set $E$ such that 
$[\varphi(\chi_{E})]=[f],$ and because $f<\varphi(\chi_{C'})$ we can take $E$ 
to be a subset of $C'.$  Set $C=C'\setminus E.$  Then $[\varphi(\chi_C)]=g$ 
and properties (1)--(3) are satisfied.
\end{proof}

\noindent{\em Proof of Theorem~\ref{sft-rp}. }Given $m\in\nats,\varepsilon>0,$ 
and a finite subset $F$ of $A$ we let $C$ be as in Lemma~\ref{towlem} and set 
$e_i=\varphi(\chi_C).$  By Remark~\ref{rem-asind}, we need not consider $F$ at 
all.  Our construction guarantees the $[e_i]$ are all the same, and that 
$[1-\sum_ie_i]<[e_0],$ so the partial isometries needed to apply 
Theorem~\ref{prohlin} exist.
\qed

\begin{proposition}If $T_1$ and $T_2$ are shift equivalent primitive integral 
matrices then\linebreak 
$A_{T_1}\times_{\alpha_{T_2}}\ints \cong A_{T_2}\times_{\alpha_{T_2}}\ints.$
\end{proposition}
\begin{proof}Let us write $\alpha_i$ for $\alpha_{T_i}$ and $A_i$ for 
$A_{T_i},$ $i=1,2.$  Shift equivalence implies that there is an isomorphism 
$\rho:A_1\cong A_2$ such that 
\[
\rho_*\circ(\alpha_1)_*=(\alpha_2)_*\circ\rho_*
\]
and
\[
\rho_*(K_0(A_1)_+)=K_0(A_2)_+.
\]
By Theorem~4.1 of \cite{EK}, there is an automorphism $\beta$ of $A_2$ 
and a unitary $u\in A_2$ such that 
\[
\alpha_1=\rho^{-1}\circ\Ad u\circ\beta\circ\alpha_2\circ\beta^{-1}\circ\rho.
\]
For $i=1,2,$ the crossed product $A_i\times_{\alpha_i}\ints$ is generated by 
$A_i$ and a unitary $v_i$ such that 
\[
v_ia=\alpha_i(a)v_i.
\]
Define $\Phi:A_1\times_{\alpha_1}\ints\to A_2\times_{\alpha_2}\ints$ by setting
\[
\Phi(v_1)=\beta^{-1}(u)v_2
\]
and
\[
\Phi(a)=\beta^{-1}(\rho(a))),\quad a\in A_1.
\]
The map $\Phi$ is a well-defined $*$-homomorphism because
\begin{align*}
\Phi(v_1a)&=\beta^{-1}(u)v_2\beta^{-1}(\rho(a)))\\
&=\beta^{-1}(u)(\alpha_2\circ\beta^{-1}\circ\rho)(a)v_2\\
&=\beta^{-1}(u(\beta\circ\alpha_2\circ\beta^{-1}\circ\rho)(a)u^*u)v_2\\
&=(\beta^{-1}\circ\rho\circ\alpha_1)(a)\beta^{-1}(u)v_2\\
&=\Phi(\alpha_1(a)v_1).
\end{align*}
To see that it is invertible observe that $\beta^{-1}\circ\rho$ is an 
isomorphism and
\[
\Phi(\rho^{-1}(u^*)v_1)
=\beta^{-1}(\rho(\rho^{-1}(u^*))\beta^{-1}(u)v_2
=\beta^{-1}(u^*u)v_2=v_2.
\]
\end{proof}
\begin{remark}
We do not know whether isomorphism of the crossed product algebras implies 
shift equivalence, or whether 
$(A_{T_1},\alpha_{T_1})\cong(A_{T_2},\alpha_{T_2})$ is equivalent to shift 
equivalence, strong shift equivalence, or neither.
\end{remark}

\section*{Acknowledgements}

This work was part of my doctoral thesis.  I thank my advisor,
Jack Wagoner, for mentorship and guidance.  I am grateful to Nate Brown, who 
introduced me to the Rohlin property for automorphisms of $C^*$-algebras.
I also thank Kenley Jung and Ian Putnam for helpful discussions.

\end{document}